\documentclass[conference]{IEEEtran}
\IEEEoverridecommandlockouts
\usepackage{cite}
\usepackage{amsmath,amssymb,amsfonts}
\usepackage{algorithmic}
\usepackage{algorithm}
\usepackage{graphicx}
\usepackage{textcomp}
\usepackage{xcolor}
\usepackage{bm}
\usepackage{hyperref}
\def\BibTeX{{\rm B\kern-.05em{\sc i\kern-.025em b}\kern-.08em
    T\kern-.1667em\lower.7ex\hbox{E}\kern-.125emX}}
\begin{document}

\title{
Optimal Unmanned Aerial Vehicle Deployment for Macro-Micro Traffic Monitoring Fused with Connected Vehicles\\
}

\author{
\IEEEauthorblockN{Chaopeng Tan}
\IEEEauthorblockA{\textit{Chair of Traffic Process Automation} \\
\textit{Technische Universität Dresden}\\
Dresden, Germany \\
chaopeng.tan@tu-dresden.de}
\and
\IEEEauthorblockN{Jiarong Yao*}
\IEEEauthorblockA{\textit{School of Electrical and Electronic Engineering} \\
\textit{Nanyang Technological University}\\
Singapore, Singapore \\
yaojoanna2018@gmail.com}
\and
\IEEEauthorblockN{Meng Wang}
\IEEEauthorblockA{\textit{Chair of Traffic Process Automation} \\
\textit{Technische Universität Dresden}\\
Dresden, Germany \\
meng.wang@tu-dresden.de}
}

\maketitle

\begin{abstract}
Reliable estimation of macro and micro traffic states is essential for urban traffic management. Unmanned Aerial Vehicles, with their airborne full-sample continuous trajectory observation, bring new opportunities for macro- and micro-traffic state estimation. In this study, we will explore the optimal UAV deployment problem in road networks in conjunction with sampled connected vehicle data to achieve more reliable estimation of macroscopic path flow as well as microscopic arrival rates and queue lengths. Oriented towards macro-micro traffic states, we propose entropy-based and area-based uncertainty measures, respectively, and transform the optimal UAV deployment problem into minimizing the uncertainty of macro-micro traffic states. A quantum genetic algorithm that integrates the thoughts of metaheuristic algorithms and quantum computation is then proposed to solve the large-scale nonlinear problem efficiently. Evaluation results on a network with 18 intersections have demonstrated that by deploying UAV detection at specific locations, the uncertainty reduction of macro-micro traffic state estimation ranges from 15.28\% to 75.69\%. A total of 5 UAVs with optimal location schemes would be sufficient to detect over 95\% of the paths in the network considering both microscopic uncertainty regarding the intersection operation efficiency and the macroscopic uncertainty regarding the route choice of road users.
\end{abstract}


\section{Introduction}
Reliable and comprehensive traffic state monitoring plays a critical role in effective urban traffic management. Conventional traffic monitoring has predominantly focused on two scales: the macro level, which involves estimating origin-destination (OD) or path flows, and the micro level, which focuses on capturing intersection dynamics such as vehicle arrivals and queue length. Macro-scale monitoring is essential for understanding the overall travel demand and network performance, yet it poses significant challenges in data collection and estimation due to the limitations of stationary sensors and the inherent difficulties in capturing complete travel patterns \cite{krishnakumari2020data, sun2024stochastic}. At the micro level, obtaining accurate data on intersection traffic flows and queue lengths is equally critical, as such data informs adaptive signal control strategies and congestion management; however, traditional approaches are often hampered by high installation costs, limited spatial coverage, and insufficient temporal resolution \cite{luo2025probabilistic, tan2021cumulative}.

Advancements in connected vehicle (CV) technology have provided a promising alternative by harnessing real-time vehicular data for traffic state estimation. Numerous studies have leveraged CV data to estimate both macro-level traffic conditions, such as OD matrices \cite{cao2021day} and path flows \cite{chen2021dynamic}, and micro-level conditions like intersection queue lengths \cite{tan2019cycle, hao2014cycle} and volumes \cite{zheng2017estimating, yao2019sampled}. Although CV data offers improved temporal resolution and spatial coverage compared to fixed sensors, they are fundamentally sampled observations and CV-based methods are thus subject to uncertainties and biases, particularly under conditions of low penetration rates. 

Simultaneously, the rapid development of unmanned aerial vehicles (UAVs) has introduced an innovative perspective to the domain of traffic monitoring. UAVs have the unique capability of continuously capturing high-resolution imagery and video over specific areas, thereby enabling the tracking of vehicle trajectories and real-time assessment of traffic conditions from a top-down view \cite{bisio2022systematic}. This aerial approach offers a distinct advantage by covering the whole area and providing a more complete picture of vehicular movements. Ideally, if we have a fleet of UAVs hovering above all intersections during a time-of-day (TOD) period for monitoring, both the macroscopic and microscopic traffic states can be perfectly obtained \cite{barmpounakis2020new}. However, such a practice is far from feasible as full coverage for all the signalized intersections may be costly regarding device procurement, configuration, and maintenance. Moreover, a 100\% coverage is sometimes redundant or unnecessary, considering the traffic flow correlation among intersections.

Given these considerations, this study uses limited UAV resources as a complement to large-scale city-wide traffic monitoring, aiming to fuse CV data and UAV data for more robust and reliable macro-micro traffic state estimation. Essentially, the research problem is to design the location of the limited UAVs to realize maximal coverage of traffic information within the network with maximal accuracy or minimal uncertainty. 
To solve this problem, a UAV location optimization framework is thus proposed to provide the exact intersection set to minimize the estimation uncertainty of traffic states at both macro- and micro-levels. The major contribution of this study is three-fold:
\begin{enumerate}
    \item As far as we know, this is the first study on collaborating CV and UAV for macro-micro traffic state monitoring. 
    \item For macroscopic and microscopic traffic states, we propose entropy-based and area-based methods to measure their uncertainty, respectively.
    \item We propose an uncertainty minimization framework that transforms the traffic state uncertainty optimization problem into a UAV location optimization problem, which improves the data utilization efficiency at the planning level and facilitates various macroscopic and microscopic traffic state estimation methods.
\end{enumerate}

\section{Problem Statement}
For a signalized urban roadway network with an intersection set $\mathcal{I}$ (indexed by $i$) and a link set $\mathcal{L}$ (indexed by $l$). To map the traffic generation point (origin) and traffic attraction point (destination) to the network, here the link pair is used to denote an OD pair, thus the OD set can be denoted as $\mathcal{OD}={(\mathcal{L_O},\mathcal{L_D} )}$ (indexed by $od$), where $\mathcal{L_O}$ and $\mathcal{L_D}$ are the original and destination links. Topologically, a path refers to the sequential link set (or intersection set with OD links) for a vehicle to travel from one origin to one destination within the network, and the path set is denoted by $\mathcal{P}$ (indexed by $p$). From the perspective of signal control, a path can also be described by a sequential movement set regarding each intersection it passes. For each intersection, the movements of left-turn, straight and right-turn are generally considered for each connecting link; here, a set $\mathcal{M}_i$ is denoted to include all the movements of intersection $i \in \mathcal{I}$, indexed by $m$. We use $\mathcal{P}_{i,m}$ to denote the path set that contains movement $m$ at intersection $i$. Then, we have
\begin{align}
    Q_{i,m} = \sum_{p \in \mathcal{P}_{i,m}} Q_p \label{eq: movement and path flow}
\end{align}
where $Q_{i,m}$ is the traffic flow of movement $m$ and $Q_p$ is the traffic flow of path $p$. This also indicates that $Q_p \leq Q_{i,m}$ for $p \in \mathcal{P}_{i,m}$.

\begin{figure}[!t]
\centering
\includegraphics[width=3.2 in]{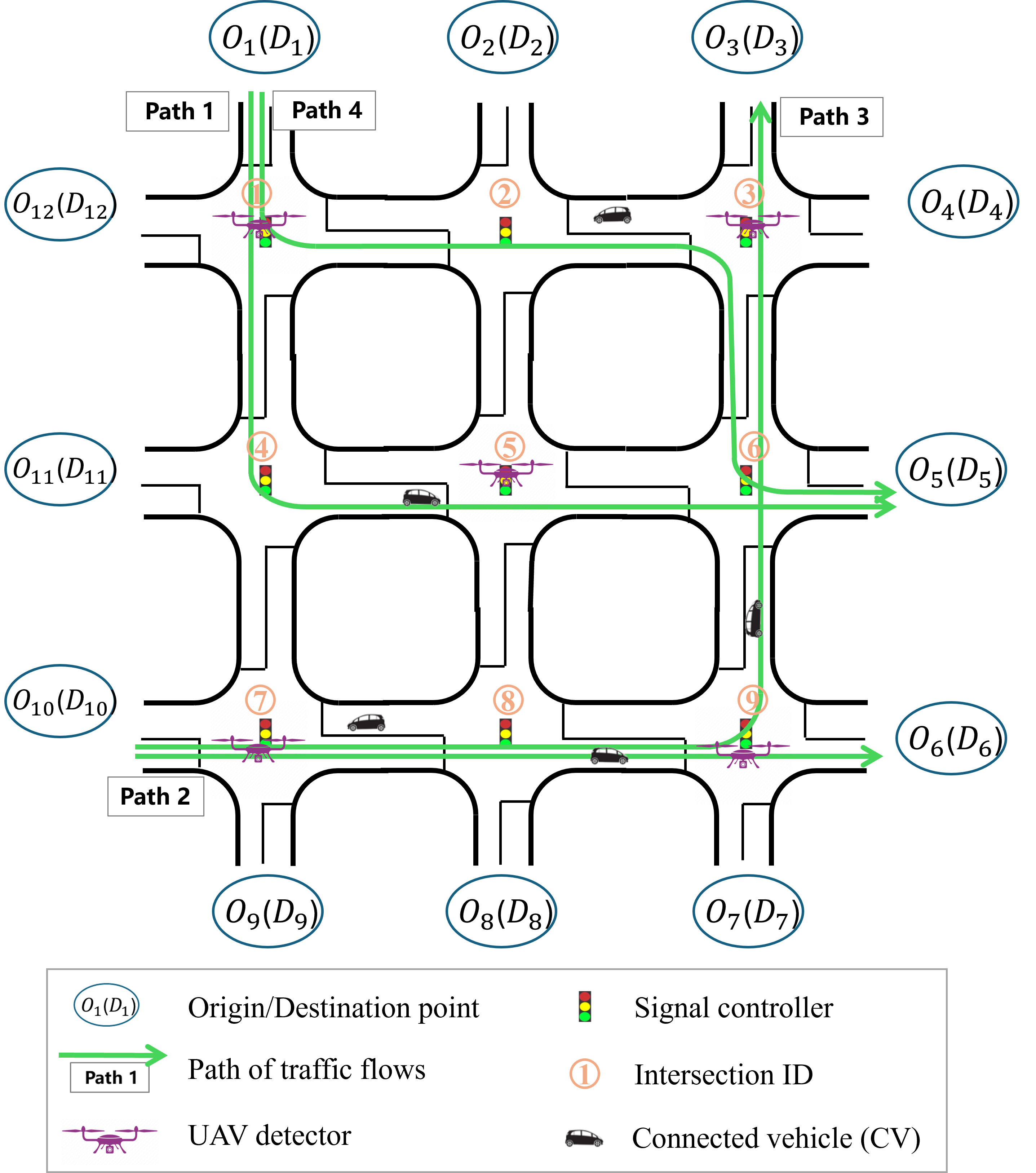}
\caption{A sketch of the research scenario.}
\label{fig: problem}
\end{figure}

Using UAV detection data at certain intersections and sampled connected vehicle trajectory data within the network scale as input, multi-resolution traffic state in a TOD period is expected to be estimated accurately, including macroscopic traffic state indicators like OD flow and path flow, and microscopic traffic state indicators such as arrival profile, queue length, control delay, number of stops etc. Let $\mathcal{I}_{uav}$ be the subset of intersections deployed with UAVs, the uncertainty of microscopic traffic state indicators of such observed intersections is zero as the UAV detector is assumed to be able to capture the whole scale of the intersection. As for intersections not covered by UAV, the microscopic traffic state indicators can only be estimated based on the sampled CVs, thus the uncertainty is closely related to the sample size and temporal distribution of CVs within each signal cycle. As for the limited number of UAVs within the network, the uncertainty of macroscopic traffic state estimation lies in to what degree the path flow and OD flow can be reconstructed from the known movement flows of the detected intersections obtained by UAVs and the sample path flow of CVs. Such uncertainty is affected by the selection of $\mathcal{I}_{uav}$ considering the travel behavior of the traffic participants in the network. Therefore, considering a signalized traffic network as shown in Fig. \ref{fig: problem}, the UAV fleet location influences the uncertainty of traffic state estimation at both macroscopic and microscopic levels. With the aim of minimizing such uncertainty, the research problem of this study is to find an optimal location scheme of the UAV fleet for reliable traffic monitoring.

Necessary assumptions are declared here:
\begin{enumerate}
    \item The number of available UAVs $N_{uav}$ is fixed and they can all be used for intersection monitoring. However, vehicles cannot be matched between different UAVs due to privacy concerns and the high altitude view of the UAVs. 
    \item The OD pair set $\mathcal{OD}$ and path set $\mathcal{P}$ are known \cite{behara2020novel}.
    \item The traffic arrivals and queue lengths at intersections where UAVs are deployed are definite, i.e., the uncertainty is 0.
    \item The historical movement flow at intersections can be estimated by historical CVs, though inaccurate \cite{zheng2017estimating}. 
\end{enumerate}

\section{Methodology}
\subsection{Location optimization model}
Given limited UAVs with prior sample CVs as auxiliaries, UAV location optimization aims to minimize the uncertainty $Z$ of macro-micro traffic state estimation of the whole network, as given by the multi-objective programming model:
\begin{align}
    \min_{I_{uav}} & Z=w_1 F_{path} + w_2 F_{queue} + w_3 F_{arrival}. \label{eq: obj}\\
    \text{s.t.} \quad & u_i = \begin{cases}
        1 \quad \text{if} \quad i \in \mathcal{I}_{uav} \\
        0 \quad \text{otherwise},
    \end{cases} \\
    & \sum_{i \in \mathcal{I}} u_i \leq N_{uav}.
\end{align}
Macroscopic uncertainty is quantified by the uncertainty of path flows $F_{path}$ as path flow represents the distribution of traffic demand within the network. Microscopic uncertainty is quantified by the uncertainties of queue length $F_{queue}$ and arrival flow $F_{arrival}$ as they demonstrate movement-level traffic demand and intersection control performance, respectively. $w_1$, $w_2$, and $w_3$ are weighting coefficients. $N_{uav}$ denotes the number of available UAVs.
Next, we introduce the calculation of $F_{path}$, $F_{queue}$, and $F_{arrival}$. 

\subsection{Entropy-based uncertainty of macroscopic traffic state}
In this section, an indicator, path flow reconstruction entropy (PFRE) is proposed to quantify the uncertainty of macroscopic path flow based on the information entropy theory. Information entropy is a measure of the uncertainty level of random variables, which describing the information needed to determine the state of the random variable on average. The bigger the information entropy, the bigger the uncertainty of the random variable. 

Recall that $Q_p$ denotes the traffic flow of path $p \in \mathcal{P}$. Considering $Q_p$ as a random variable, the set of possible values of $Q_p$ is denoted as $\mathcal{Q}_p$. Then, the information entropy of path $p$ is given by
\begin{align}
    H(p) = -\sum_{q \in\mathcal{Q}_p} Pr(q)\log_2 Pr(q),
\end{align}
where $Pr(q)$ denotes the probability of $Q_p = q$.
Given UAV deployment $\mathcal{I}_{uav}$, the conditional entropy of path $p$ is
\begin{align}
    H(p \mid \mathcal{I}_{uav}) = -\sum_{q \in\mathcal{Q}_{p\mid uav}} Pr(q \mid \mathcal{I}_{uav})\log_2 Pr(q \mid \mathcal{I}_{uav})
\end{align}
where $Pr(q \mid \mathcal{I}_{uav})$ is the conditional probability of $Q_p = q$ given $\mathcal{I}_{uav}$ and the corresponding set of possible values becomes $\mathcal{Q}_{p \mid uav}$. The total information entropy $H(\mathcal{P} \mid \mathcal{I}_{uav})$ of all paths within the network is 
\begin{align}
    H(\mathcal{P} \mid \mathcal{I}_{uav}) = \sum_{p \in \mathcal{P}} H(p \mid \mathcal{I}_{uav}).
\end{align}

Next, we derive the conditional entropy $Pr(q \mid \mathcal{I}_{uav})$. We use $\mathcal{I}_p$ to denote the set of intersections on path $p$. Then we have set $\mathcal{I}_{uav,p} = \mathcal{I}_{uav} \cap \mathcal{I}_{p}$ and $\mathcal{I}_{uav,p}^- = \mathcal{I}_p - \mathcal{I}_{uav,p}$, where $\mathcal{I}_{uav,p}$ denotes intersections on path $p$ deployed with UAVs and $\mathcal{I}_{uav,p}^-$ denotes the other intersections on path $p$ without UAVs. 

Regarding intersections on path $p$, if the intersection $i \in \mathcal{I}_{uav,p}$, then we can use both CV and UAV observations to determine a tighter interval of the path flow $Q_p$, thus reducing its uncertainty; otherwise, only CV observations can provide information.  

We use $m_p \in \mathcal{M}_i$ to denote the movement index of path $p$ at intersection $i$. The corresponding number of CVs is denoted as $N_{i,m_p}$. The number of CVs of path $p$ is denoted as $N_{p}$. 
Similar to \eqref{eq: movement and path flow}, we have $N_{i,m_p} = \sum_{p' \in \mathcal{P}_{i,m_p}} N_{p'}$. 
The proportion $r_{m_p}$ of path $p$ to movement flow of $m_p$ is obtained as
\begin{align}
    r_{m_p} = N_{p}/N_{i,m_p}.
\end{align}
Considering flow uncertainty, we can use historical data from multiple periods to obtain 95th percentile confidence intervals for $r_{m_p}$, i.e., $[r_{m_p}^l, r_{m_p}^u]$. 

For intersection $i \in \mathcal{I}_{uav,p}$, we can get an interval of the corresponding path flow of the path $p$ as
\begin{align}
    [Q_{p,i}^l, Q_{p,i}^u] = [N_{p}, Q_{m_p} r_{m_p}^u],
\end{align}
where $Q_{m_p}$ is the historical movement flow of $m_p \in \mathcal{M}_i$ estimated by historical CVs. This can be interpreted as the interval that is expected to be obtained after deploying a UAV.
While for intersection $i \in \mathcal{I}_{uav,p}^-$, we only have CV information for the path flow interval of path $p$ as
\begin{align}
    [Q_{p,i}^l, Q_{p,i}^u] = [N_{p}, \infty].
\end{align}
Then, regarding path $p$, the set of possible values $\mathcal{Q}_{p\mid uav}$ of its path flow given $\mathcal{I}_{uav,p}$ is obtained as 
\begin{align}\label{eq:path flow feasible zone}
    \mathcal{Q}_{p\mid uav} &= \left[ \max_{i \in \mathcal{I}_p} \{Q_{p,i}^l\}, \min_{i \in \mathcal{I}_p} \{Q_{p,i}^u\} \right]  \nonumber \\ 
    &= \left[N_{p},  \min_{i \in \mathcal{I}_{uav,p}} \{Q_{m_p} r_{m_p}^u\}\right],
\end{align}
where only integers are considered. Since we have no information about the distribution of the values in $\mathcal{Q}_{p\mid uav}$, we assume that the values follow a uniform distribution. The conditional probability $Pr(q \mid \mathcal{I}_{uav})$ for all $q \in \mathcal{Q}_{p\mid uav}$ is calculated as
\begin{align}
    Pr(q \mid \mathcal{I}_{uav}) = 1/\left|\mathcal{Q}_{p\mid uav}\right|,
\end{align}
where $|\cdot|$ denotes the size of the set. 

To summarize, the PFRE, i.e., $H(\mathcal{P} \mid \mathcal{I}_{uav})$, of the network is calculated as
\begin{align}
    F_{path} &= H(\mathcal{P} \mid \mathcal{I}_{uav}) \nonumber \\
    &= \sum_{p \in \mathcal{P}} \sum_{q \in\mathcal{Q}_{p\mid uav}} \frac{1}{\left|\mathcal{Q}_{p\mid uav}\right|}\log_2 \left|\mathcal{Q}_{p\mid uav}\right|.
\end{align}

\subsection{Area-based uncertainty of microscopic traffic state}
In this section, an indicator, uncertainty area (UA), is proposed to quantify the uncertainty of the microscopic traffic states at intersections, namely, queue length and traffic arrivals. Considering the shockwave evolution externalized by the motion status change of CVs, the uncertain range of cycle-based queue length and arrival rate estimation can be reduced by CV trajectories, regardless of the specific theoretical estimation models. 

\subsubsection{Uncertainty of Queue length} 
\begin{figure}[!t]
\centering
\includegraphics[width=3.2 in]{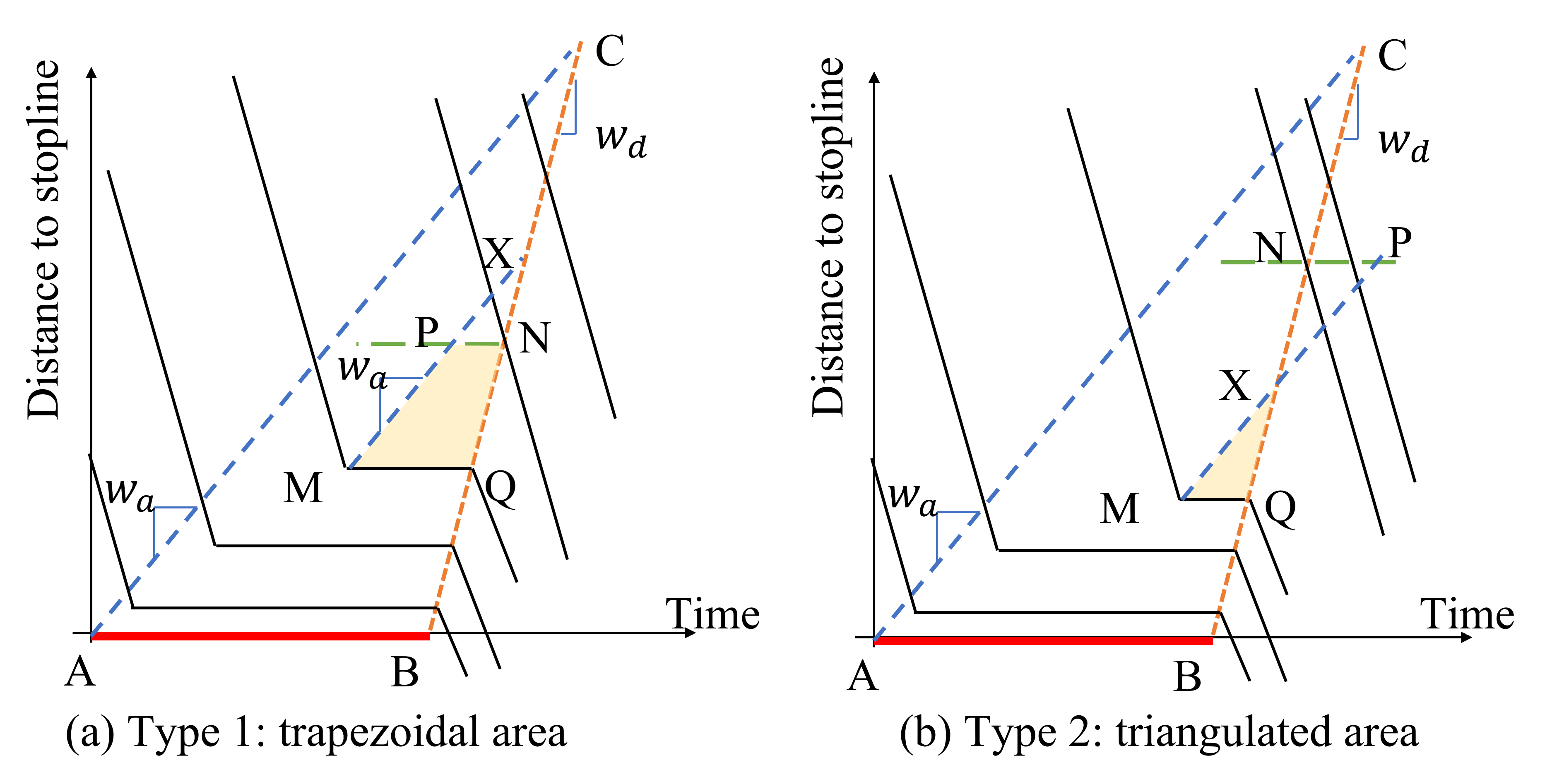}
\caption{Uncertain area of queue length.}
\label{fig: uncertainty of queue}
\end{figure}

For a specific cycle $k$ of the movement $m$ at intersection $i$, the uncertain area of its queue length is shown in Fig. \ref{fig: uncertainty of queue}. For brevity, indices $k, m, i$ are omitted. Given the signal timing parameters of cycles, the maximum queue accumulation wave speed $w_a$, and the queue dissipation wave speed $w_d$, we can define a global feasible area $S_{global,q}$ of possible queue length, i.e., the triangle area $ABC$: 
\begin{align}
    S_{global,q} = 0.5 w_a w_d R^2/(w_d-w_a)
\end{align}
where $R$ is the red time of the cycle.

Basically, the final feasible queue length is determined by two points: the joining queue point $M$ of the last queued CVs and the crossing point $N$ of the line segment $\overline{BC}$ and the first non-queued CV’s trajectory (no matter whether this CV succeeds in passing in the current cycle). The coordinates of point $M$ are denoted by $M:(t_M,d_M )$, while the coordinates of point $N$ are denoted by $N:(t_N,d_N)$, both of which can be easily solved based on simple plane geometry. Note that for the unsaturated and oversaturated traffic conditions, the differences in the model derivation lies in the process of obtaining $M$ and $N$. The subsequent steps are identical.

Given point $M$, we can generate a line parallel to $\overline{AC}$ and intersecting with $\overline{BC}$ at point $X:(t_X,d_X )$, and we can solve that
\begin{align}
    &t_X=(w_d R-w_a t_M+d_M)/(w_d-w_a ), \\
    &d_X=(w_d w_a (R-t_X )+d_X)/(w_d-w_a ).
\end{align}

Then, based on the relative positions of $X$ and $N$, the uncertain area of queue length has two types of shapes: 
\begin{itemize}
    \item \textbf{Type 1: trapezoidal area.} When $d_X>d_N$, the uncertain area of queue length, i.e., $S_{u,q}$, is a trapezoid. Point $Q:(t_Q,d_Q )$ is also known based on the trajectory of the last queued CV. Point $P:(t_P,d_P)$ is solved as $t_P=t_M+(d_N-d_M)/w_a$ and $d_P=d_N$. Then we have 
    \begin{align}
        S_{u,q}=0.5 (t_Q-t_M+t_N-t_P )(d_N-d_Q )
    \end{align}
    \item \textbf{Type 2: triangulated area.} When $d_X \leq d_N$, $S_{u,q}$ is a triangle. We have \begin{align}
        S_{u,q}=0.5(t_Q-t_M )(d_X-d_Q ).
    \end{align}
\end{itemize}

The uncertainty area refers to the probability of traffic state estimates given sampled CV data, which is denoted by $U_{queue}$ for the queue length:
\begin{align}
    U_{queue} = S_{u,q}/S_{global,q}.
\end{align}

Easily, we can find that $U_{queue}<1$ always exists and when a UAV is deployed at the intersection, $U_{queue} = 0$. Thus, we have $0 \leq U_{queue} < 1$. 
In particular, the physical meaning of $U_{queue}$ is the uncertainty regarding the spatial-temporal plane based on the observed CV trajectories, even before any estimation model is applied. Therefore, this index is independent of any theoretical model for queue length estimation and is a direct evaluation of CV data when applied to queue length estimation. A smaller $U_{queue}$ indicates a smaller possible area of queue length, suggesting that any method is expected to achieve a more accurate queue length estimate based on the observed CV trajectories.

\subsubsection{Uncertainty of vehicle arrivals}
Vehicle arrivals are generally quantified as the arrival profile, and the arrival rate indicates the number of arrival vehicles per second in a cycle; thus, the global feasible area of arrival rate $S_{global,a}$ should be a rectangle in the arrival rate-time plane:
\begin{align}
    S_{global,\lambda} = \lambda_u C,
\end{align}
where $\lambda_u$ is the maximum arrival rate and $C$ is the cycle length. 

\begin{figure}[!t]
\centering
\includegraphics[width=3.4 in]{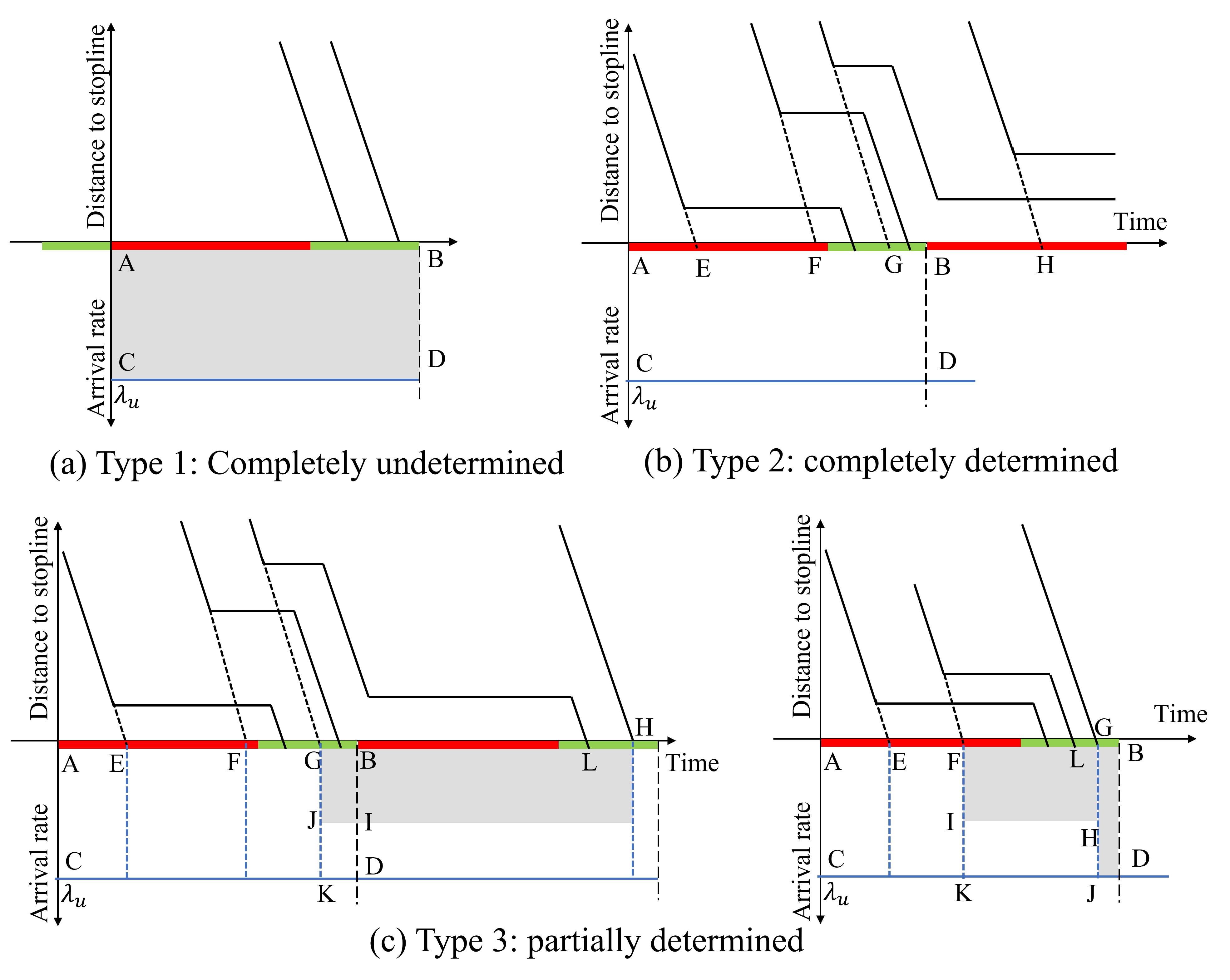}
\caption{Uncertain area of vehicle arrivals.}
\label{fig: uncertainty of arrival}
\end{figure}

The queuing process of CVs at intersections can provide partial vehicle arrival information during the cycle. For example, by using the queuing positions and expected arrival time of two consecutive queued CVs, we can easily estimate the arrivals of non-CVs in between. Therefore, there are three types of uncertain areas of arrival rates given CV observations during the cycle, as shown in Fig. \ref{fig: uncertainty of arrival}. 
\begin{itemize}
    \item \textbf{Type 1: Completely undetermined.} If only non-queued CVs or no CVs are observed during the cycle, then the arrival rate of this cycle is completely indeterminate. The uncertain area is 
    \begin{align}
        S_{u,\lambda} = S_{global,\lambda}. 
    \end{align}
    \item \textbf{Type 2: Completely determined.} If a cycle has twice-queued CVs (who queued for the first time in this cycle) and the next cycle has queued CVs, then the arrival rate during this cycle is completely determined, since vehicle arrivals between all queued CVs are known. We have
    \begin{align}
        S_{u,\lambda} = 0.
    \end{align}
    \item \textbf{Type 3: Partially determined.} Except for type 1 and 2 cycles, all the rest of the cycles belong to type 3, whose arrival rate during the cycle is partially determined. Generally, there are two sub-types: a) with twice-queued CVs in the cycle and only non-queued CVs captured in the next cycle, and b) with non-queued CVs in the cycle.    
    For the first sub-type, arrival rates during $T_{AE}$, $T_{EF}$, and $T_{FG}$ are known, while arrival rates during $T_{GH}$ have an upper bound $\lambda_{GH,ub}=T_{HL}/(T_{GH} h_s )$. Then, the uncertain area of arrival rate of the is 
    \begin{align}
        S_{u,\lambda} = \lambda_{GH,ub} (C-t_G).
    \end{align}
    For the second sub-type, arrival rates $T_{AE}$ and $T_{EF}$ are known, arrival rates during $T_{FG}$ has an upper bound $\lambda_{FG,ub}=T_{GL}/(T_{FG} h_s )$, and arrival rates during $T_{GB}$ is completely unknown, thus we have
    \begin{align}
        S_{u,\lambda} =\lambda_u (C-t_G )+\lambda_{FG,ub} (t_G-t_F )
    \end{align}
\end{itemize}

Similar to the $U_{queue}$, we can also define the uncertainty index for arrival rate, i.e., $U_{arrival}$, as below
\begin{align}
    U_{arrival} = S_{u,\lambda}/S_{global,\lambda}
\end{align}
We can find that $0 \leq U_{arrival} \leq 1$. Note that, $U_{arrival}=0$ does not mean that estimates by existing methods are error-free. Multi-lane scenarios may introduce errors due to adjacent trajectories not belonging to the same lane. A smaller $U_{arrival}$ indicates less uncertainty of the arrival profile, suggesting that any method is expected to achieve a more accurate arrival rate estimate based on the observed CV trajectories.

In summary, given $\mathcal{I}_{uav}$ for UAV deployment, the microscopic uncertainty of the whole network is calculated as
\begin{align}
    F_{queue} = \sum_{i \in \mathcal{I}}  u_i \sum_{m \in \mathcal{M}_i} \sum_{k \in \mathcal{K}_i} U_{queue}^{i,m,k} \\
    F_{arrival} = \sum_{i \in \mathcal{I}}  u_i \sum_{m \in \mathcal{M}_i} \sum_{k \in \mathcal{K}_i} U_{arrival}^{i,m,k} 
\end{align}
where $\mathcal{K}_i$ denotes the set of all cycles of intersection $i$ during the analysis period. 

\subsection{Solution algorithm}
As a variant of the set coverage problem, the proposed UAV location optimization model is actually an NP-hard problem, which is hard to solve using exact algorithms within acceptable computation time, thus prompting the usage of metaheuristic algorithms for an efficient and effective solution. Especially for the research scenario in this study, with the increase of the roadway network scale, the problem presents nonlinearity with a large number of local extremes, which calls for better diversity representation and evolution randomness. A quantum genetic algorithm (QGA) \cite{malossini2008quantum} is adopted here to obtain the solution integrating the thoughts of metaheuristic algorithms and quantum computation. 

As the core of quantum-based metaheuristic algorithm, quantum computation actually provides a new computation mode using qubit for information processing in place of binary digits. Qubit is the basic information unit in quantum coding, represented by the superposition state of two basis states, $|0\rangle$ and $|1\rangle$, as given by Eq. \eqref{eq: qubit}, showing that a qubit can be an arbitrary value between $|0\rangle$ and $|1\rangle$.

\begin{align}
    |\Psi\rangle = \alpha |0\rangle + \beta |1\rangle \label{eq: qubit} \\
    \alpha^2 + \beta^2 = 1
\end{align}

$\alpha, \beta$ are complex numbers, denoting the probability amplitudes of the basis states. $\alpha^2,\beta^2$ represent the probability of being collapsed to the state of $|0\rangle$ or $|1\rangle$, respectively. 

Accordingly, the element chromosome of the population is represented by qubits, and the quantum state collapse of all the qubits is for population evolution according to the changes of individual fitness. A chromosome containing $I$ qubits can be described by Eq. \eqref{eq: chromosome}. Each qubit $q_i \quad (i=1,2,\dots, I)$ of $Chrom$ will collapse into one certain state for each observation of the chromosome. For the proposed model, the state of $q_i$ represents whether a UAV is placed over intersection $i$ for detection, i.e., $q_i$ with state $|1\rangle$ means intersection $i$ is chosen to set up UAV detection, otherwise, $q_i$ with state $|0\rangle$ means UAV detection is not available at intersection $i$.  

\begin{align}
    Chrom =  \left[
\begin{vmatrix} \alpha_{1} \\ \beta_{1} \end{vmatrix}
\begin{vmatrix} \alpha_{2} \\ \beta_{2} \end{vmatrix}
\dots
\begin{vmatrix} \alpha_{I} \\ \beta_{I} \end{vmatrix}
\right] \label{eq: chromosome} 
\end{align}

As the objective of the proposed method is to minimize the macro-micro traffic state estimation uncertainty, here $-Z$ is used as the fitness function. The pseudocode is shown in Algorithm \ref{alg:AQGA}.

\begin{algorithm}
\caption{QGA for UAV Location}
\label{alg:AQGA}
\begin{algorithmic}[1]
\REQUIRE Network topology $(\mathcal{I},\mathcal{L})$, sample CV data, UAV fleet size $N_{uav}$, chromosome population size $P_c$, chromosome length $I_c$, number of evolution generations $T$
\ENSURE Optimal UAV location $\mathcal{I}_{uav}$
\STATE Input parameters $N_{p}, \{U_{queue}^{i,m,k}\}, \{U_{arrival}^{i,m,k}\}$
\STATE Generate initial population: $P(0) = \{\text{Chrom}_1^{(0)}, \text{Chrom}_2^{(0)}, \dots, \text{Chrom}_{M_c}^{(0)}\}$
\STATE Initialize each chromosome individual at generation $t=0$. The chromosome length $I_c$ is exactly the number of intersections of the network, and each gene represents whether the specific intersection is detected by UAV.
\STATE Measure the observation value of the initial population and obtain observation values $\vec{u}(t,p) = \{u_1^{t,p}, u_2^{t,p}, \dots, u_i^{t,p}, \dots\}$
\FORALL{$\vec{u}(t,p)$}
    \STATE Calculate the fitness value $\text{fit}(\vec{u}(t,p))$
\ENDFOR
\STATE $\vec{u}^{t,\text{best}} = \arg\max_{p \in P_c} \text{fit}(\vec{u}(t,p))$
\WHILE{$t < T$}
    \STATE Evolve $t$ to $t+1$ using quantum rotation gate\cite{malossini2008quantum}: $t = t + 1$
    \STATE Repeat lines 5 to 7
    \STATE $B(t) = \arg\max_{p \in P_c} \text{fit}(\vec{u}(t,p))$
    \IF{$\text{fit}(\vec{u}^{t,\text{best}}) < \text{fit}(B(t))$}
        \STATE $\vec{u}^{t,\text{best}} = B(t)$
    \ENDIF
\ENDWHILE
\end{algorithmic}
\end{algorithm}

\section{Evaluation}
The proposed method is evaluated through a simulation case based on an empirical network in Qingdao, China to justify both the performance of multi-resolution traffic state uncertainty evaluation and solution algorithm.  
\subsection{Study site}
\begin{figure}[!t]
\centering
\includegraphics[width=3.2 in]{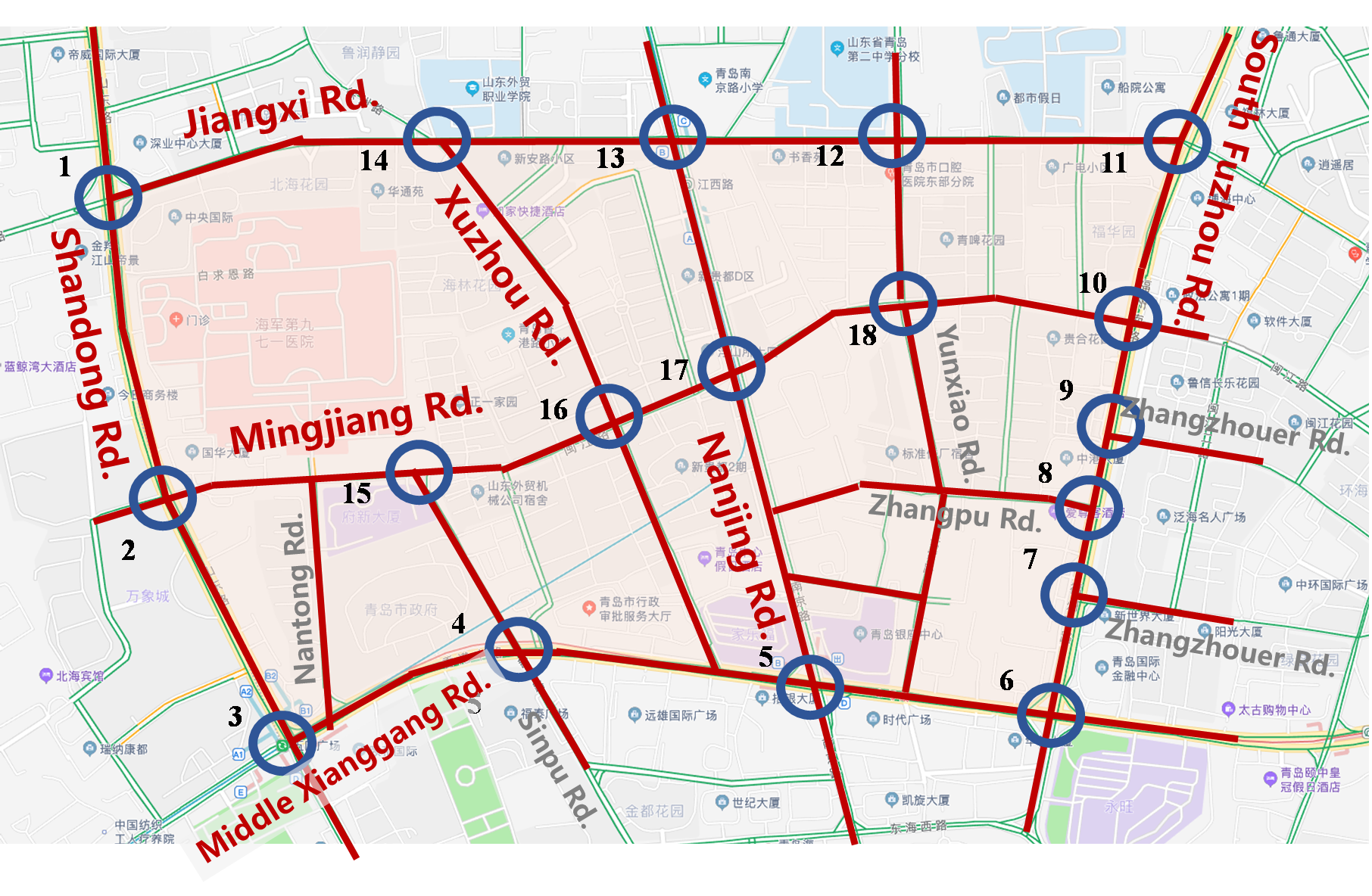}
\caption{Study site.}
\label{fig: shinan_network}
\end{figure}

The case study site is chosen as the "4 Longitudinal and 3 Transverse" roadway network in Shinan District, Qingdao, Shandong Province. As shown in Fig. \ref{fig: shinan_network}, there are 18 signalized intersections (marked out by blue circles with index in black), including 14 typical crossroads and 4 T-type intersections. Through VISSIM software, the simulation model of this Shinan network was built and the empirical data from 7:00 – 9:00 on Mar. 1st, 2019, including license plate recognition (LPR) data, floating car data (FCD) and signal timing data were used for calibration. Based on empirical detection data, there were in total 28 origins (destinations) and 311 paths. 

The calibrated simulation model was run for three times, with a simulation period of 2 hours each time, to extract the vehicle trajectories at an uploading frequency of 1s. A penetration rate of 0.1 was applied to obtain the sampled trajectories as the historical CVs to estimate the movement flows \cite{yao2019sampled}, uncertainty area of queue as well as arrival rate at each intersection, and the average values of three runs were used as the historical movement flow $Q_{m_p}$,  $U_{queue}^{i,m,k}$ and $ U_{arrival}^{i,m,k}$. 

\subsection{Uncertainty evaluation}
Using the CV data with a penetration rate of 0.1 and the given number of available UAVs as input, the proposed method was tested with $N_{uav}$ ranging from 1 to 18. Coded in Python, the method was run on a laptop with an AMD Ryzen 7 5800H CPU @ 3.20GHz, 16.0 GB of RAM, and one NVIDIA GeForce RTX 3060 GPU. For the QGA solver, the numbers of population and generation were set as 20 and 200 respectively. As for the weight coefficient, we set $w_1 = w_2 = w_3$. 

\begin{figure}[!t]
\centering
\includegraphics[width=3.4 in]{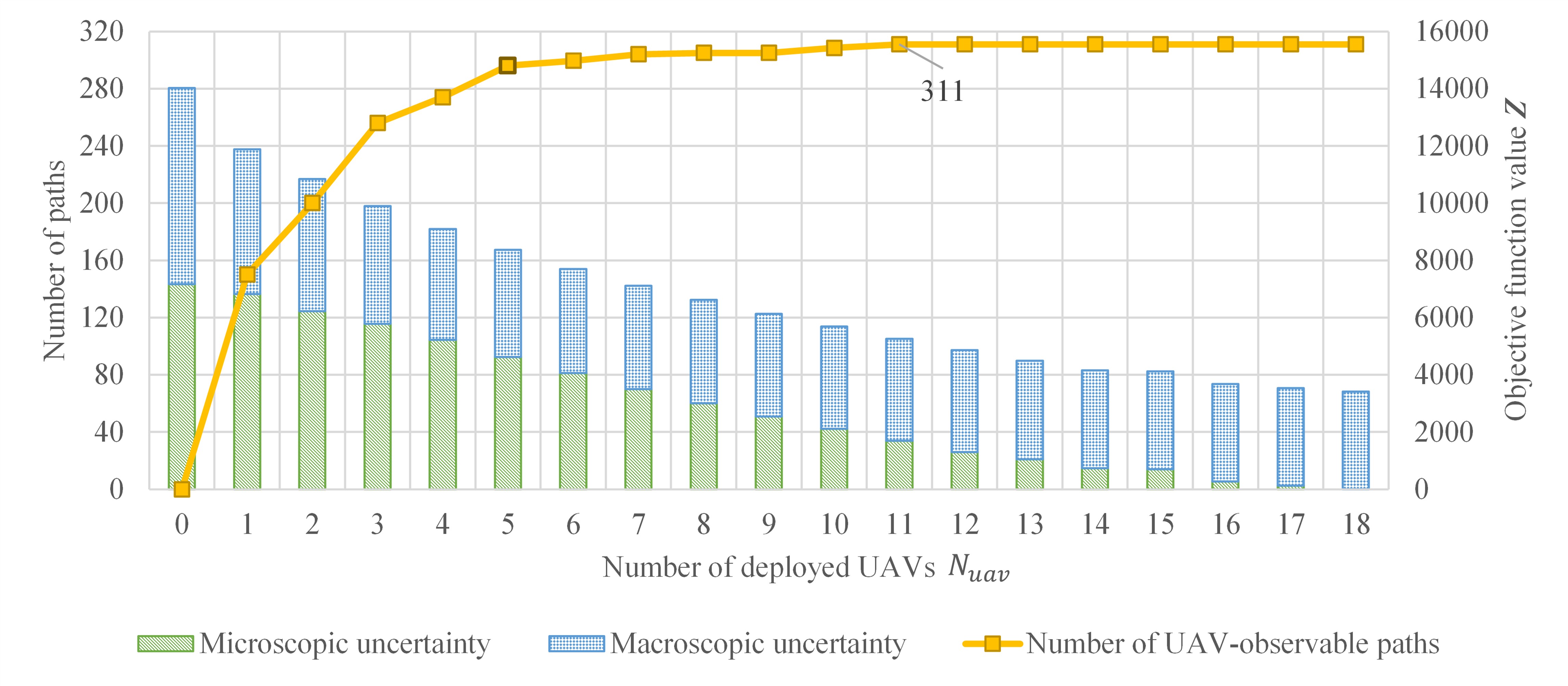}
\caption{Uncertainty using different UAV numbers.}
\label{fig: sensitivity of UAV num}
\end{figure}

As shown in Fig. \ref{fig: sensitivity of UAV num}, the uncertainty of macro-micro traffic state estimation of the whole network does decrease with the deployment of the UAV number, by about 75.69\% from 0 UAV to full coverage, and the decreasing amplitude shows a decreasing trend from 15.28\% to 3.48\%. As for the uncertainty of the two sub-objectives, the decrease amplitude of microscopic uncertainty grows with the increase of UAV-detected intersection, while the decrease of macroscopic uncertainty slowly drops with the increase of UAV number, showing a more evident marginal reduction as the UAV number is less than 6. Even when all the intersections in the whole network are covered by UAVs, the uncertainty of traffic state estimation still exists,  although only in the macroscopic side. It is sensible that even the movement flows of all the intersections are known, the estimation of path flow still requires the information of route choice of drivers which is partially known through the sampled CV data. As long as the flow value set of each path of Eq. \eqref{eq:path flow feasible zone} is an interval, the macroscopic uncertainty of path flow reconstruction still exists. 

The yellow line in Fig. \ref{fig: sensitivity of UAV num} represents the time-varying trend of the number of paths detected by UAVs. Similar to the uncertainty index, the number of UAV-observable paths also shows a growing trend when the number of UAV-detected intersections increases. All 311 paths can be observable when the number of UAV-detected intersections increases to 11 while the marginal increase starts to flatten when the UAV number is larger than 5, which also accrods with the trend of the above-mentioned macroscopic uncertainty. Thus, 5 UAVs can be regarded as a cost-effective choice of the UAV device number to be deployed in the studied network, as more than 95\% of the network paths can be observable with only less than one-third of the intersections detected by UAVs. 

\begin{figure}[!t]
\centering
\includegraphics[width=3.4 in]{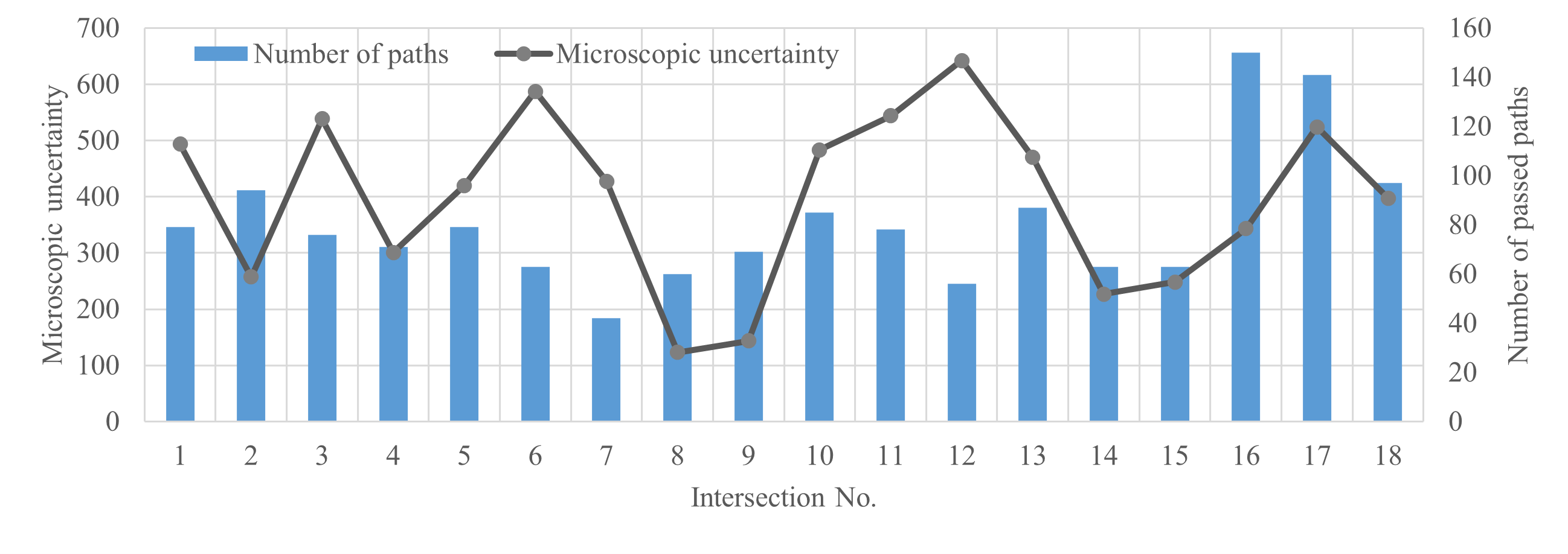}
\caption{Information of intersections.}
\label{fig: intersection information}
\end{figure}

The optimal UAV location scheme is to assign the 5 UAVs to Intersection 1, 3, 6, 12, and 16 as shown in Fig. \ref{fig: shinan_network}. Considering the microscopic uncertainty and the number of paths passing each intersection shown in Fig. \ref{fig: intersection information}, both the intersections with large uncertainty and the intersection with the most paths passing through are chosen to be deployed with UAV detection. Based on the equal weights and the relative magnitudes of macroscopic and microscopic uncertainty shown in Fig. \ref{fig: sensitivity of UAV num}, this corresponding optimal UAV location scheme also demonstrates that the solution algorithm effectively works to search the desired multi-objective optimum.

\begin{figure}[!t]
\centering
\includegraphics[width=3 in]{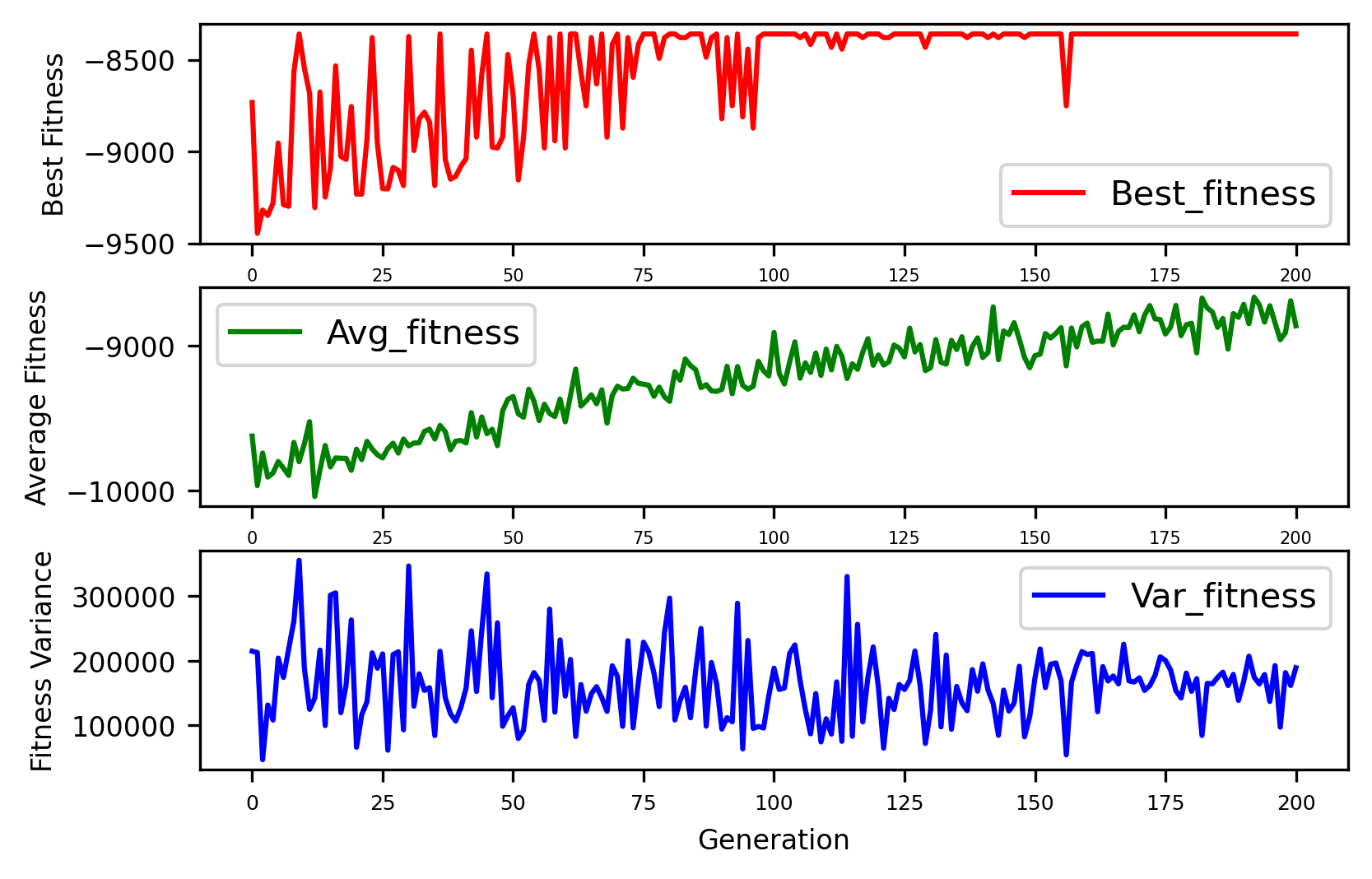}
\caption{Algorithm performance of 5 UAVs.}
\label{fig: result of 5 UAVs}
\end{figure}

Fig. \ref{fig: result of 5 UAVs} further shows the algorithm performance of QGA to obtain the optimal solution. It can be seen from the upper red line chart of the time-varying trend of the best fitness of each population that although the algorithm obtained the best solution at the 9th generation at first, the fitness value still fluctuated severely. With the evolution of the population, the algorithm started to converge after the 100th generation and finally converged at the 156th generation. As the opposite number of $Z$ is used as the fitness function, the increasing trend of the average fitness also implies that the population is getting better for exploitation and exploration of the optimum. With the population evolves, the variance also becomes smaller and less fluctuate. 

\begin{figure}[!t]
\centering
\includegraphics[width=3 in]{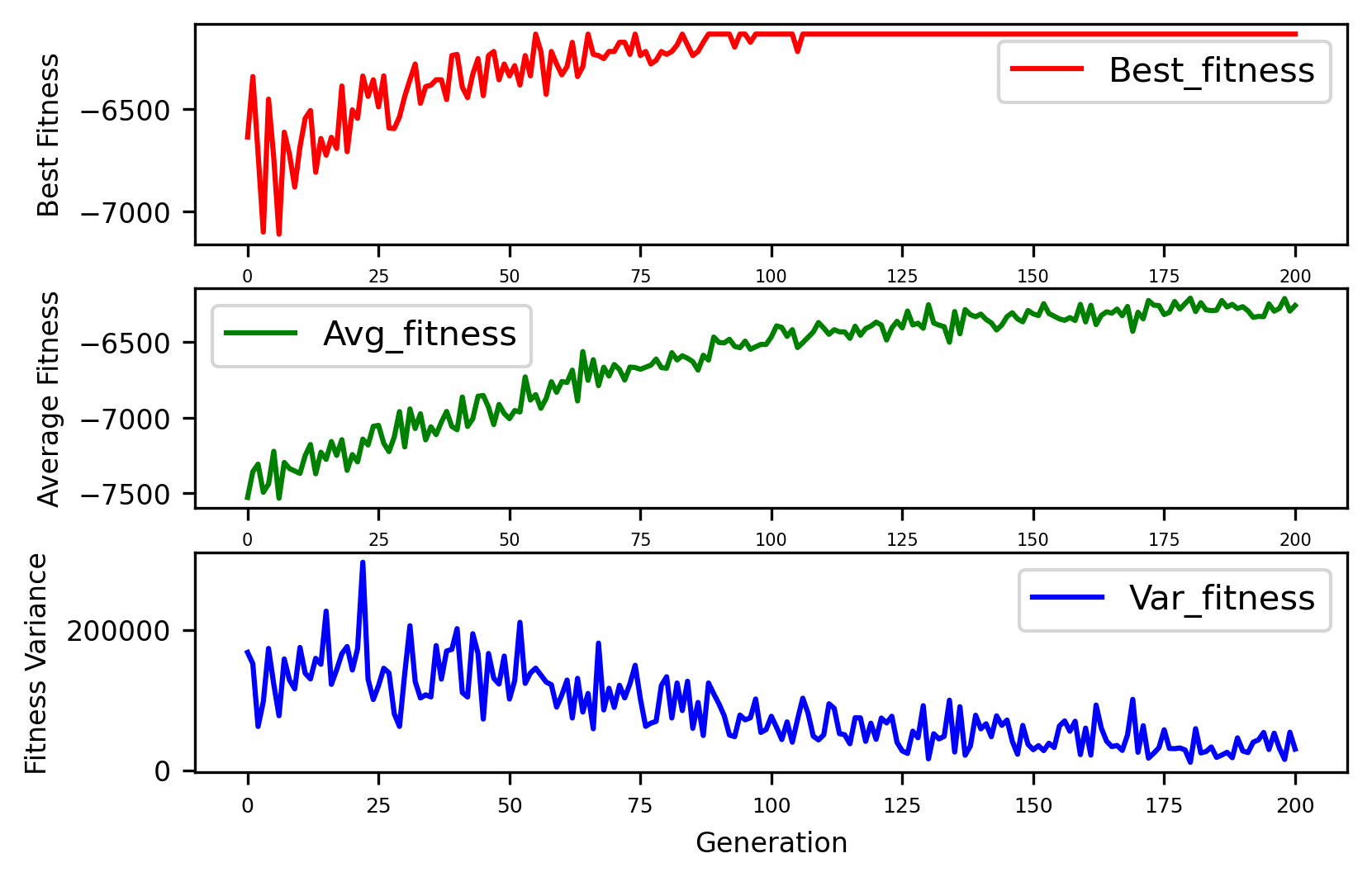}
\caption{Algorithm performance of 9 UAVs.}
\label{fig: result of 9 UAVs}
\end{figure}

With the increase of available UAVs to be deployed, the performance of QGA also shows slight difference in population variance and convergence. The case of 5-UAV and 13-UAV share the same size of the feasible zone, the performance of fitness variance is similar in Fig. \ref{fig: result of 5 UAVs} and Fig. \ref{fig: result of 5 UAVs}. As for the case of 9-UAV which has the largest size of the feasible zone for the study site, the results in Fig. \ref{fig: result of 9 UAVs} shows that QGA converges faster with smaller population-wise fluctuation. Besides, the decrease of population variance of 9-UAV case is more obvious than those of 5-UAV case and 13-UAV case, which implies that QGA may focus more on exploration of optimum searching in small solution space while more on exploitation in large solution space.

\begin{figure}[!t]
\centering
\includegraphics[width=3 in]{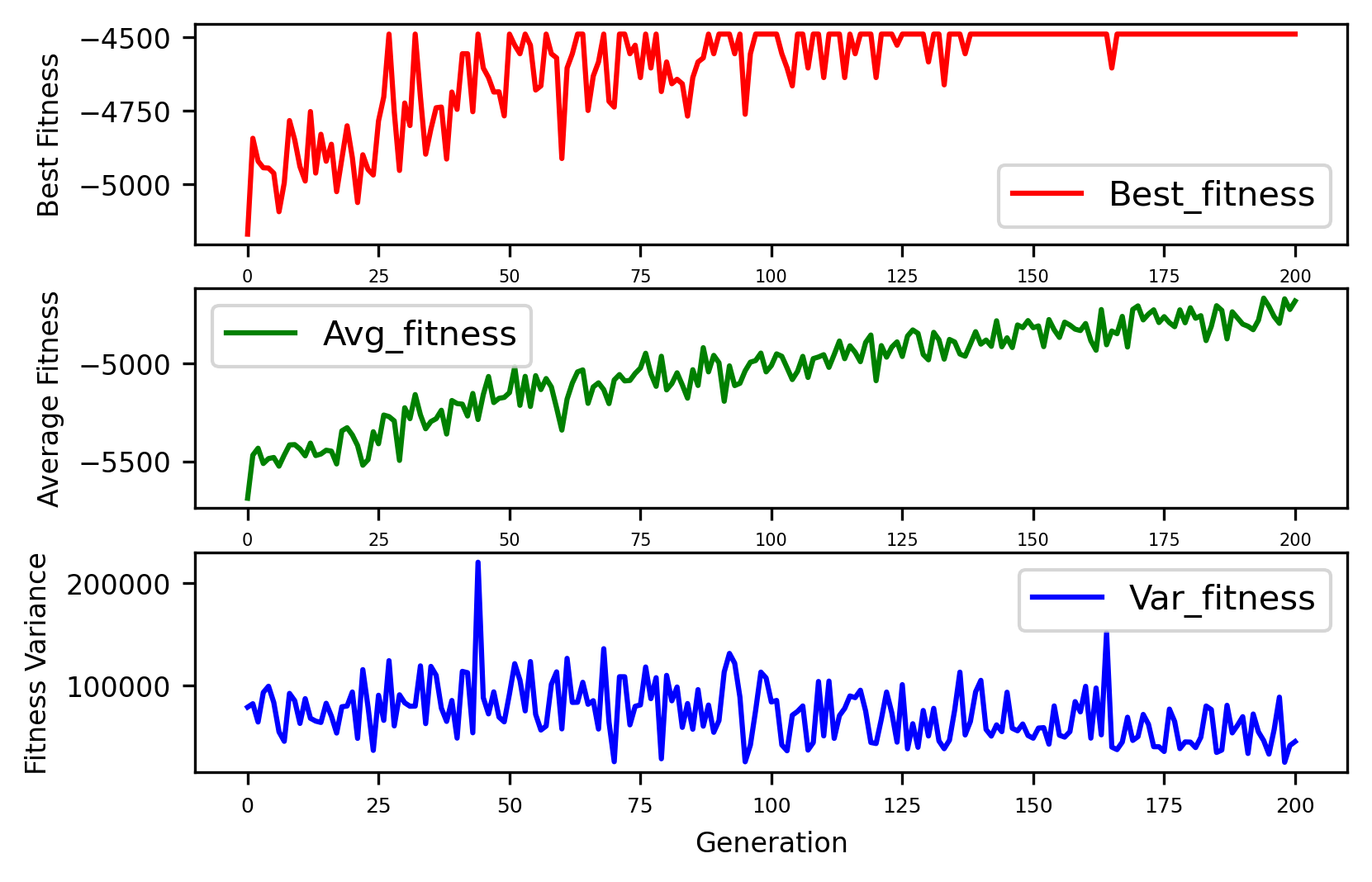}
\caption{Algorithm performance of 13 UAVs.}
\label{fig: result of 13 UAVs}
\end{figure}

\section{Conclusion and Future Work}
This study explores the optimal UAV deployment problem in road networks in conjunction with sampled connected vehicle data to achieve a more reliable estimation of macroscopic path flow as well as microscopic arrival rates and queue lengths. Specifically, the optimal UAV deployment problem is transformed into an optimization problem to minimize the uncertainty of macro-micro traffic states, where entropy-based and area-based approaches are proposed to measure the uncertainty of macro-micro traffic states, respectively. Evaluation results have demonstrated that deploying a small number of UAVs at specific locations can significantly reduce the uncertainty of macro-micro traffic states. A total of 5 UAVs with optimal location schemes would be sufficient to detect over 95\% of the paths in the network and reduce the network uncertainty by about 40.3\%, significantly beneficial for both macro-micro traffic state estimation. 

Future research includes i) sensitivity analysis on factors like the weights taken in Eq. \eqref{eq: obj}, CV penetration rate as well as parameters of QGA and ii) further testing of the enhancement of UAV deployment on different macro and micro traffic state estimation methods.

\bibliographystyle{IEEEtran}
\bibliography{Manuscript_1006}{}

\end{document}